# $L_2$ boosting in kernel regression

B.U. PARK[1], Y.K. LEE[2] and S. HA[1]

[1]*Department of Statistics, Seoul National University, Seoul 151-747, Korea.*
[2]*Department of Statistics, Kangwon National University, Chuncheon 200-701, Korea.*

In this paper, we investigate the theoretical and empirical properties of $L_2$ boosting with kernel regression estimates as weak learners. We show that each step of $L_2$ boosting reduces the bias of the estimate by two orders of magnitude, while it does not deteriorate the order of the variance. We illustrate the theoretical findings by some simulated examples. Also, we demonstrate that $L_2$ boosting is superior to the use of higher-order kernels, which is a well-known method of reducing the bias of the kernel estimate.

## 1. Introduction

In the last decade, several important approaches for classification and pattern recognition have been proposed with feasible computational algorithms in the machine learning community. Boosting is one of the most promising techniques that has recently received a great deal of attention from the statistical community. It was first proposed by Schapire (1990) as a means of improving the performance of a given method, called a *weak learner*. Subsequent investigations of the methods have been made in both communities. These include, among others, Freund (1995); Freund and Schapire (1996, 1997); Schapire, Freund, Bartlett and Lee (1998); Breiman (1998, 1999); Schapire and Singer (1999); Friedman, Hastie and Tibshirani (2000); Friedman (2001).

Understanding boosting algorithms as functional gradient descent techniques gives theoretical justifications of the methods; see Mason, Baxter, Bartlett and Frean (2000) and Friedman (2001). It connects various boosting algorithms to statistical optimization problems with corresponding loss functions. For example, AdaBoost (Freund and Schapire (1996)) can be interpreted as giving an approximate solution, starting from an initial learner, to the problem of minimizing the exponential risk for classification. Also, Logit-Boost corresponds to an approximate optimization of the log-likelihood of binary random variables, see Friedman *et al.* (2000).

In this paper, we study boosting as a successful bias reduction method in nonparametric regression. Since the regression function $m(\cdot) = E(Y|X = \cdot)$ is the minimizer of the $L_2$ risk $E[Y - m(X)]^2$, it is natural to take the squared error loss as an objective function.

Application of the functional gradient descent approach to the $L_2$ risk is trivial since minimization of the $L_2$ risk itself is a linear problem and thus there is no need to linearize it. In fact, the population version of $L_2$ boost is nothing else than adding $m - m_0$ to an initial function $m_0$ so that a single update yields an exact solution. However, the empirical $L_2$ boost is non-trivial. It amounts to repeated least-squares fitting of residuals.

$L_2$ boost in the context of regression has been studied by Friedman (2001) and Bühlmann and Yu (2003). In the latter work, the authors provided some expressions for the average squared bias and the average variance of the $L_2$ boost estimate obtained from a linear smoother in terms of the eigenvalues of the corresponding smoother matrix. They also showed that, if the learner is a smoothing spline, it is possible for $L_2$ boosting to achieve the optimal rate of convergence for all higher-order smoothness of the regression function. In doing so, they took the iteration number, rather than the penalty constant, as the regularization parameter. The optimal rate is attained if one takes the iteration number $r = O(n^{2p/(2\nu+1)})$ as the sample size $n$ goes to infinity, where $p$ is the order of the smoothing spline learner and $\nu$ is the smoothness of the regression function.

In this paper, we investigate the theoretical and empirical properties of $L_2$ boosting when the learner is the Nadaraya–Watson kernel smoother. We derive the bias and variance properties of the estimate in terms of the bandwidth (smoothing parameter), which is more conventional in nonparametric function estimation. We show that the optimal rate of convergence is also achieved by the Nadaraya–Watson $L_2$ boosting for all smoothness of the regression function if the iteration number $r$ is high enough, depending on the smoothness $\nu$, and the bandwidth is properly chosen as $O(n^{-1/(2\nu+1)})$. In particular, we prove that each step of $L_2$ boosting reduces the bias of the estimate by two orders of the bandwidth, and also that additional boosting steps do not deteriorate the order of the variance. We illustrate these theoretical findings by some simulated examples in a numerical study. Also, we compare the finite sample properties of $L_2$ boosting with those of higher-order kernel smoothing, the latter being a well-known method of reducing the bias of the estimate. Our results suggest that $L_2$ boosting is superior to the use of higher-order kernels.

## 2. Main results

The $L_2$ boosting algorithm is derived from application of the functional gradient descent technique to the $L_2$ loss. The task of the latter is to find the function $m$ that minimizes a functional $\psi(m)$. With an initial function $m_0$, one searches the best direction $\delta$ such that $\psi(m_0 + \varepsilon\delta)$ is minimized. Let $\dot{\psi}(\delta)$ be the Gâteaux differential of $\psi$ with increment $\delta$, that is,

$$\dot{\psi}(\delta) = \lim_{\varepsilon \to 0} \frac{\psi(\cdot + \varepsilon\delta) - \psi(\cdot)}{\varepsilon}.$$

To first order in $\varepsilon$, minimizing $\psi(m_0 + \varepsilon\delta)$ with respect to $\delta$ is equivalent to minimizing $\dot\psi(\delta)(m_0)$. Let $\delta_1$ denote the minimizer. The update of the initial $m_0$ is given by $m_1 = m_0 + \varepsilon_1\delta_1$, where $\varepsilon_1$ minimizes $\psi(m_0 + \varepsilon\delta_1)$. Then, the process is iterated.

Let $m(x) = E(Y|X = x)$ be the regression function. If one applies the functional gradient descent technique to the $L_2$ loss, $\psi(m) = \frac{1}{2}E[Y - m(X)]^2$, then one gets

$$\dot\psi(\delta)(m_0) = -E[\delta(X)(Y - m_0(X))].$$

Since minimizing $-E[\delta(X)(Y - m_0(X))]$ subject to $E\delta(X)^2 = c$ for some constant $c > 0$ is equivalent to minimizing $E[Y - m_0(X) - \delta(X)]^2$, it follows that the updated function $m_1$ is given by $m_1 = m_0 + \delta_1$, where

$$\delta_1 = \underset{\delta}{\operatorname{argmin}}\, E[Y - m_0(X) - \delta(X)]^2 = E[Y - m_0(X)|X = \cdot].$$

Thus, with the $L_2$ loss, the update $m_1$ equals the true function $m$:

$$m_1(x) = m_0(x) + E[Y - m_0(X)|X = x] = m(x).$$

The $L_2$ boosting algorithm given below is an empirical version of the updating procedure above.

### *Algorithm ($L_2$ Boosting).*

Step 1 *(Initialization): Given a sample $\mathcal{S} = \{(X_i, Y_i), i = 1,\ldots,n\}$, fit an initial estimate $\hat m_0(x) \equiv \hat m(x; \mathcal{S})$ to the data.*
Step 2 *(Iteration): Repeat for $r = 1,\ldots,R$.*

  (i) *Compute the residuals $e_i = Y_i - \hat m_{r-1}(X_i)$, $i = 1,\ldots,n$.*
  (ii) *Fit an estimate $\hat m(x; \mathcal{S}_e)$ to the data $\mathcal{S}_e = \{(X_i, e_i), i = 1,\ldots,n\}$.*
  (iii) *Update $\hat m_r(x) = \hat m_{r-1}(x) + \hat m(x; \mathcal{S}_e)$.*

Thus, $L_2$ boosting is simply repeated least-squares fitting of residuals. With $r = 1$ (one-step boosting), it has been already proposed by Tukey (1977), usually referred to as "twicing". Twicing is related to using higher-order kernels. It was observed by Stützle and Mittal (1979) that, in the case of the fixed equispaced design points $x_i = i/n$, twicing a kernel smoother is asymptotically equivalent to directly using a higher-order kernel. To be more specific, let $K$ be a kernel function, $h > 0$ be the bandwidth and $K_h(u) = K(u/h)/h$. Define $K^* = 2K - (K * K)$, where $*$ denotes the convolution operator. Note that $K^*$ is a higher-order kernel. If $\hat m_0(x) = n^{-1}\sum_{i=1}^n K_h(x - x_i)Y_i$, then $\hat m_1(x) \simeq n^{-1}\sum_{i=1}^n K_h^*(x - x_i)Y_i$, where $\simeq$ is due to the integral approximation error $n^{-1}\sum_{j=1}^n K_h(x - x_j)K_h(x_j - x_i) \simeq \int K_h(x - z)K_h(z - x_i)\,\mathrm{d}z = (K * K)_h(x - x_i)$.

In this paper, we consider random covariates $X_i$. We derive the theoretical properties of $L_2$ boosting when the learner is the Nadaraya–Watson kernel smoother, that is,

$$\hat m(x; \mathcal{S}) = \frac{\sum_{i=1}^n K_h(X_i - x)Y_i}{\sum_{i=1}^n K_h(X_i - x)}.$$

The Nadaraya–Watson smoothing is the simplest and numerically most stable technique of local kernel regression. We note that statistical properties of $L_2$ boosting for $r \geq 1$ with Nadaraya–Watson smoothing have not been investigated before.

Throughout the paper, we assume $K$ is a symmetric probability density function which is Lipschitz continuous and is supported on $[-1, 1]$. The bounded support condition for $K$ can be relaxed to include kernels, such as Gaussian, that decrease to zero at tails with an exponential rate. Below, we discuss the asymptotic properties of $\hat{m}_r$ for $r \geq 1$. For this, we assume that $h \to 0$ and $nh/\log n \to \infty$ as $n \to \infty$.

We denote a pre-estimate of $m$ by $\widetilde{m}$. Thus, at the $r$th iteration $\widetilde{m} = \hat{m}_{r-1}$. Let $\hat{m}$ be its update defined by

$$\hat{m}(x) = \widetilde{m}(x) + \sum_{i=1}^{n} w_i(x)[Y_i - \widetilde{m}(X_i)], \tag{1}$$

where $w_i(x) = [\sum_{j=1}^{n} K_h(X_j - x)]^{-1} K_h(X_i - x)$. At the $r$th iteration $\hat{m} = \hat{m}_r$. Note that, for the initial estimate $\hat{m}_0(x) = \sum_{j=1}^{n} w_j(x) Y_j$, we get

$$\hat{m}_0(x) - m(x) = \sum_{j=1}^{n} w_j(x) e_j + \sum_{j=1}^{n} w_j(x)[m(X_j) - m(x)],$$

where $e_j = Y_j - m(X_j)$.

Let $\widetilde{w}_j$ be the weight functions for $\widetilde{m}$ that depend solely on $X_1, \ldots, X_n$ and satisfy

$$\sum_{j=1}^{n} \widetilde{w}_j(x) = 1 \quad \text{for all } x, \qquad \widetilde{m}(x) = \sum_{j=1}^{n} \widetilde{w}_j(x) Y_j.$$

Define the updated weight functions by

$$\hat{w}_j(x) = w_j(x) + \widetilde{w}_j(x) - \sum_{i=1}^{n} w_i(x) \widetilde{w}_j(X_i).$$

Note that $\hat{w}_j$ also depends solely on $X_1, \ldots, X_n$. One can verify

$$\sum_{j=1}^{n} \hat{w}_j(x) = 1 \quad \text{for all } x, \qquad \hat{m}(x) = \sum_{j=1}^{n} \hat{w}_j(x) Y_j,$$

so that

$$\hat{m}(x) - m(x) = \sum_{j=1}^{n} \hat{w}_j(x) e_j + \sum_{j=1}^{n} \hat{w}_j(x)[m(X_j) - m(x)]. \tag{2}$$

From (2), we note that $\mathrm{Var}(\hat{m}(x)|X_1, \ldots, X_n) = \sum_{j=1}^{n} \hat{w}_j(x)^2 \sigma^2(X_j)$, where $\sigma^2(x) = \mathrm{Var}(Y|X = x)$. The following theorem provides the magnitude of the conditional variance. Let $f$ denote the marginal density $f$ of the covariate $X_i$. We assume that $f$ is supported on $\mathcal{I} = [0, 1]$.

**Theorem 1.** *Assume that $f$ is continuous on $\mathcal{I}$ and $\inf_{x \in \mathcal{I}} f(x) > 0$. If $\sup_{x \in \mathcal{I}} \sum_{j=1}^{n} \widetilde{w}_j(x)^2 = O_p(n^{-1}h^{-1})$ uniformly for $x \in \mathcal{I}$, then $\sup_{x \in \mathcal{I}} \sum_{j=1}^{n} \hat{w}_j(x)^2 = O_p(n^{-1}h^{-1})$ uniformly for $x \in \mathcal{I}$.*

In the proof of Theorem 1 given below, we prove that $\sum_{j=1}^{n} w_j^2(x) = O_p(n^{-1}h^{-1})$ uniformly for $x \in \mathcal{I}$. Thus, the weight functions $w_j$ for the initial estimate $\hat{m}_0$ satisfy the condition of Theorem 1. This means that $L_2$ boosting does not deteriorate the order of the variance of the estimate as the iteration goes on.

**Proof of Theorem 1.** It follows that

$$\sum_{j=1}^{n} \hat{w}_j^2(x) \leq 3 \sum_{j=1}^{n} w_j^2(x) + 3 \sum_{j=1}^{n} \widetilde{w}_j^2(x) + 3 \sum_{j=1}^{n} \left[ \sum_{i=1}^{n} w_i(x) \widetilde{w}_j(X_i) \right]^2$$

$$\leq 3 \sum_{j=1}^{n} w_j^2(x) + 3 \sum_{j=1}^{n} \widetilde{w}_j^2(x) + 3 \sum_{j=1}^{n} \sum_{i=1}^{n} w_i(x) \widetilde{w}_j(X_i)^2$$

$$\leq 3 \sum_{j=1}^{n} w_j^2(x) + 3 \left[ \sup_{x \in \mathcal{I}} \sum_{j=1}^{n} \widetilde{w}_j(x)^2 \right] \left[ 1 + \sum_{i=1}^{n} w_i(x) \right]$$

$$= 3 \sum_{j=1}^{n} w_j^2(x) + 6 \left[ \sup_{x \in \mathcal{I}} \sum_{j=1}^{n} \widetilde{w}_j(x)^2 \right].$$

To complete the proof, it remains to show that $\sum_{j=1}^{n} w_j^2(x) = O_p(n^{-1}h^{-1})$ uniformly for $x \in \mathcal{I}$. Let $\mathcal{I}_h = [h, 1-h]$. Then,

$$n^{-1} \sum_{i=1}^{n} K_h(X_i - x) = \begin{cases} f(x) + o_p(1), & \text{uniformly for } x \in \mathcal{I}_h, \\ f(x) C_1(x) + o_p(1), & \text{uniformly for } x \in \mathcal{I}/\mathcal{I}_h, \end{cases}$$

$$n^{-1} h \sum_{i=1}^{n} [K_h(X_i - x)]^2 = \begin{cases} f(x) C_2 + o_p(1), & \text{uniformly for } x \in \mathcal{I}_h, \\ f(x) C_3(x) + o_p(1), & \text{uniformly for } x \in \mathcal{I}/\mathcal{I}_h, \end{cases}$$

where $1/2 \leq C_1(x) \leq 1$, $C_2 = \int_{-1}^{1} K^2$ and $\int_{0}^{1} K^2 \leq C_3(x) \leq \int_{-1}^{1} K^2$. From this, we conclude

$$\sum_{j=1}^{n} w_j^2(x) = \frac{1}{nh} \frac{n^{-1} h \sum_{j=1}^{n} [K_h(X_i - x)]^2}{[n^{-1} \sum_{i=1}^{n} K_h(X_i - x)]^2} = O_p(n^{-1}h^{-1})$$

uniformly for $x \in \mathcal{I}$. □

Next, we discuss the conditional bias of the update $\hat{m}$. The conditional biases of $\widetilde{m}$ and $\hat{m}$ equal $\sum_{j=1}^{n} \widetilde{w}_j(x)[m(X_j) - m(x)]$ and $\sum_{j=1}^{n} \hat{w}_j(x)[m(X_j) - m(x)]$, respectively.

In the case where the pre-estimate $\widetilde{m}$ is the initial estimate $\hat{m}_0$, we have $\widetilde{w}_j = w_j$ and

$$\sum_{j=1}^{n} w_j(x)[m(X_j) - m(x)] = \frac{EK_h(X_1 - x)[m(X_1) - m(x)]}{EK_h(X_1 - x)} + \mathrm{O}_p\left(\sqrt{\frac{h \log n}{n}}\right)$$

uniformly for $x \in \mathcal{I}_\varepsilon$ for arbitrarily small $\varepsilon > 0$, sufficient smoothness of $m$ and $f$ permitting. The $\mathrm{O}_p(\sqrt{n^{-1} h \log n})$ in the above expansion comes from the mean zero stochastic terms in the numerator and denominator of the left-hand side.

**Theorem 2.** *Assume that $f$ is continuously differentiable on $\mathcal{I}_- = (0,1)$ and $\inf_{x \in \mathcal{I}} f(x) > 0$. Let $r \geq 1$ be an integer. Suppose that*

$$\sum_{j=1}^{n} \widetilde{w}_j(x)[m(X_j) - m(x)] = h^{2r} \alpha_n(x) + \mathrm{O}_p\left(\sqrt{\frac{h \log n}{n}}\right) \tag{3}$$

*uniformly for $x \in \mathcal{I}_\varepsilon$ for arbitrarily small $\varepsilon > 0$, where $\alpha_n$ is a sequence of functions that are twice differentiable on $\mathcal{I}_-$ and satisfies*

$$\lim_{\delta \to 0} \limsup_{n \to \infty} \sup_{|u-v| \leq \delta} |\alpha_n^{(k)}(u) - \alpha_n^{(k)}(v)| = 0 \tag{4}$$

*for $k = 0, 1, 2$. Then,*

$$\sum_{j=1}^{n} \hat{w}_j(x)[m(X_j) - m(x)] = h^{2(r+1)} \beta_n(x) + \mathrm{O}_p\left(\sqrt{\frac{h \log n}{n}}\right)$$

*uniformly for $x \in \mathcal{I}_\varepsilon$ for arbitrarily small $\varepsilon > 0$, where $\beta_n(x)$ is a deterministic sequence such that*

$$\beta_n(x) = -\frac{1}{2}\left[\frac{\alpha_n''(x) f(x) + 2 \alpha_n'(x) f'(x)}{f(x)}\right] \int u^2 K(u) \, du + \mathrm{o}(1).$$

Theorem 2 tells that each step of $L_2$ boosting improves the asymptotic bias of the estimate by two orders of magnitude if $m$ and $f$ are sufficiently smooth. When $\widetilde{m} = \hat{m}_0$,

$$\begin{aligned}\alpha_n(x) &= h^{-2} \frac{EK_h(X_1 - x)[m(X_1) - m(x)]}{EK_h(X_1 - x)} \\ &= \frac{1}{2}\left[\frac{m''(x) f(x) + 2 m'(x) f'(x)}{f(x)}\right] \int u^2 K(u) \, du + \mathrm{o}(1),\end{aligned} \tag{5}$$

which can be shown to satisfy (4), sufficient smoothness of $m$ and $f$ permitting. In general, if $m$ and $f$ are sufficiently smooth, the corresponding sequence of the functions $\alpha_n$ at each step of the iteration satisfies (4).

For the functions class

$$\mathcal{F}(\nu, C) = \{m : |m^{(\lfloor \nu \rfloor)}(x) - m^{(\lfloor \nu \rfloor)}(x')| \leq C|x - x'|^{\nu - \lfloor \nu \rfloor} \text{ for all } x, x' \in \mathcal{I}\},$$

where $\lfloor \nu \rfloor$ is the largest integer that is less than $\nu$, it is known that the minimax optimal rate of convergence for estimating $m$ equals $n^{-\nu/(2\nu+1)}$. Let $\hat{m}_r$ denote the estimate updated at the $r$th iteration. The following theorem implies that the $L_2$ boosted Nadaraya–Watson estimate is minimax optimal if the iteration number $r$ is high enough and the bandwidth is chosen appropriately.

**Theorem 3.** *Assume that* $m \in \mathcal{F}(\nu, C_1)$, $f \in \mathcal{F}(\nu - 1, C_2)$ *for* $\nu \geq 2$ *and* $\inf_{x \in \mathcal{I}} f(x) > 0$. *Let* $r \geq \lfloor \nu/2 \rfloor$ *be an integer. Then,*

$$E[\hat{m}_r(x)|X_1, \ldots, X_n] - m(x) = \mathrm{O}_p\left(h^\nu + \sqrt{\frac{h \log n}{n}}\right) \quad (6)$$

*uniformly for* $x \in \mathcal{I}_\varepsilon$ *for arbitrarily small* $\varepsilon > 0$.

Theorems 1 and 3 imply that

$$E[(\hat{m}_r(x) - m(x))^2 | X_1, \ldots, X_n] = \mathrm{O}_p(n^{-1}h^{-1} + h^{2\nu})$$

for $r \geq \lfloor \nu/2 \rfloor$. Thus, if one takes $h = \mathrm{O}(n^{-1/(2\nu+1)})$ and $r \geq \lfloor \nu/2 \rfloor$, then $\hat{m}_r$ achieves the minimax optimal rate of convergence. We note that Bühlmann and Yu (2003) obtained similar results for smoothing spline learners. They took the iteration number as the regularization parameter and held the penalty constant fixed. In the case of the cubic smoothing spline learner, for example, they showed that if $r = \mathrm{O}(n^{4/(2\nu+1)})$, then the $r$th updated estimate achieves the optimal rate; see their Theorem 3.

**Proof of Theorem 2.** Fix $\varepsilon > 0$. Then, for sufficiently large $n$, all $X_i$ with $\sup_{x \in \mathcal{I}_\varepsilon} w_i(x) > 0$ lie in $\mathcal{I}_{\varepsilon/2}$. Thus, the expansion (3) holds if we replace $x$ by a random $X_i$ with $\sup_{x \in \mathcal{I}_\varepsilon} w_i(x) > 0$. This implies that, uniformly for $x \in \mathcal{I}_\varepsilon$,

$$\sum_{i=1}^n w_i(x) \left[\sum_{j=1}^n \widetilde{w}_j(X_i)(m(X_j) - m(X_i))\right] = h^{2r} \sum_{i=1}^n w_i(x)\alpha_n(X_i) + \mathrm{O}_p(\rho_n), \quad (7)$$

where $\rho_n = \sqrt{n^{-1}h \log n}$. From (3) and (7), we have

$$\sum_{j=1}^n \hat{w}_j(x)[m(X_j) - m(x)] = \sum_{j=1}^n \widetilde{w}_j(x)(m(X_j) - m(x))$$

$$- \sum_{i=1}^n w_i(x)\left[\sum_{j=1}^n \widetilde{w}_j(X_i)(m(X_j) - m(X_i))\right] \quad (8)$$

$$= h^{2r} \sum_{i=1}^{n} w_i(x)[\alpha_n(x) - \alpha_n(X_i)] + \mathrm{O}_p(\rho_n)$$

uniformly for $x \in \mathcal{I}_\varepsilon$. Define $\gamma_n(u,x) = [\alpha_n(x) - \alpha_n(u)]f(u)$. Then,

$$\frac{1}{n} \sum_{i=1}^{n} K_h(X_i - x)[\alpha_n(x) - \alpha_n(X_i)] = \int K_h(u-x)\gamma_n(u,x)\,\mathrm{d}u + \mathrm{O}_p(\rho_n)$$

$$= \frac{1}{2} h^2 \gamma_n''(x,x) \int u^2 K(u)\,\mathrm{d}u + r_n(x) + \mathrm{O}_p(\rho_n)$$

uniformly for $x \in \mathcal{I}_\varepsilon$, where

$$r_n(x) = h^2 \int_{-1}^{1} \int_0^1 u^2 K(u)[\gamma_n''(x - huv, x) - \gamma_n''(x,x)](1-v)\,\mathrm{d}v\,\mathrm{d}u$$

and $\gamma_n''(u,x) = \partial^2 \gamma_n(u,x)/\partial u^2$. Note that $\gamma_n''(x,x) = -[\alpha_n''(x)f(x) + 2\alpha_n'(x)f'(x)]$ and that, for any $\delta > 0$,

$$\limsup_{n \to \infty} \sup_{x \in \mathcal{I}_\varepsilon} h^{-2}|r_n(x)| \leq \limsup_{n \to \infty} \sup_{x \in \mathcal{I}_\varepsilon} \sup_{|u| \leq \delta} |\gamma_n''(x-u,x) - \gamma_n''(x,x)| \int u^2 K(u)\,\mathrm{d}u.$$

Thus, from (4) we obtain $r_n(x) = \mathrm{o}(h^2)$ uniformly for $x \in \mathcal{I}_\varepsilon$. Since $n^{-1} \sum_{i=1}^{n} K_h(X_i - x) = f(x) + \mathrm{o}_p(1)$ uniformly for $x \in \mathcal{I}_\varepsilon$, we complete the proof of Theorem 2. □

**Proof of Theorem 3.** Let $p = \lfloor \nu/2 \rfloor$. When $p = 0$ ($\nu = 2$), we know

$$E[\hat{m}_0(x)|X_1, \ldots, X_n] - m(x) = h^2 \alpha_n(x) + \mathrm{O}_p(\rho_n),$$

where $\alpha_n$ is given at (5). When $p \geq 1$ ($\nu > 2$), one can verify by repeated applications of Theorem 2 that

$$E[\hat{m}_{p-1}(x)|X_1, \ldots, X_n] - m(x) = h^{2p} \alpha_n(x) + \mathrm{O}_p(\rho_n)$$

uniformly for $x \in \mathcal{I}_\varepsilon$ for arbitrarily small $\varepsilon > 0$, where $\alpha_n$ is a sequence of functions. If $\nu = 2p + \xi$ for some integer $p \geq 1$ and $0 < \xi \leq 1$, then $\alpha_n$ satisfies

$$\limsup_{n \to \infty} \sup_{|u-v| \leq \delta} |\alpha_n(u) - \alpha_n(v)| \leq C_1 \delta^\xi$$

for some $C_1 > 0$. Since

$$E[\hat{m}_p(x)|X_1, \ldots, X_n] - m(x) = h^{2p} \sum_{i=1}^{n} w_i(x)[\alpha_n(x) - \alpha_n(X_i)] + \mathrm{O}_p(\rho_n) \qquad (9)$$

uniformly for $x \in \mathcal{I}_\varepsilon$ as in (8), we obtain (6).

Next, if $\nu = 2p + 1 + \xi$ for some integer $p \geq 1$ and $0 < \xi \leq 1$, then

$$\limsup_{n \to \infty} \sup_{|u-v| \leq \delta} |\alpha'_n(u) - \alpha'_n(v)| \leq C_2 \delta^\xi$$

for some $C_2 > 0$. Note that

$$\left| n^{-1} \sum_{i=1}^{n} K_h(X_i - x)[\alpha_n(x) + \alpha'_n(x)(X_i - x) - \alpha_n(X_i)] \right|$$

$$\leq C_2 h^\xi n^{-1} \sum_{i=1}^{n} |X_i - x| K_h(X_i - x) \tag{10}$$

$$= O_p(h^{1+\xi})$$

uniformly for $x \in \mathcal{I}_\varepsilon$. From (9) and (10), we obtain (6) in the case $\nu = 2p + 1 + \xi$, too. This completes the proof of Theorem 3. □

Two important issues that need particular attention are the choice of the bandwidth $h$ and that of the iteration number $r$, which may have substantial influence on the performance of the estimator for a finite sample size. These are related to each other in the sense that both $h$ and $r$ are regularization parameters and interplay each other. An optimal choice for one of them depends on the choice of the other. In their smoothing spline approach, Bühlmann and Yu (2003) fixed the penalty constant, whose role is the same as that of the bandwidth $h$ in our setting, and find the optimal rate of increase for $r$ (as the sample size grows), as given in the above paragraph. Our theory is for the other way around. It suggests that taking sufficiently large $r$ so that $r \geq \lfloor \nu/2 \rfloor$, but fixed without tending to infinity as the sample size grows, gives an optimal performance in terms of rate of convergence if the bandwidth $h$ is chosen in an optimal way.

In practical implementation of the boosting algorithm where the sample size is fixed, letting $r \to \infty$ alone leads to overfitting and thus jeopardizes the boosting method. One may think it is possible to avoid overfitting by increasing the bandwidth. However, increasing the bandwidth to reduce the variance of the estimator would also increase the bias, which may result in an increase of the mean squared error if $r$ is too high. Thus, one should use a data-dependent stopping rule for the iteration, as well as a data-driven bandwidth selector. For this one may employ a cross-validatory criterion, or the test bed method, as discussed in Györfi, Kohler, Krzyzak and Walk (2002) and Bickel, Ritov and Zakai (2006). To describe the latter method for selection of both $h$ and $r$, write $\hat{m}_r(\cdot;h)$ rather than $\hat{m}_r$ to stress its dependence on $h$, and let $\{(X_{n+1}, Y_{n+1}), \ldots, (X_{n+B}, Y_{n+B})\}$ be a test bed sample that is independent of the training sample $\{(X_1, Y_1), \ldots, (X_n, Y_n)\}$. Define, for each $r \geq 1$,

$$\hat{h}_r = \operatorname{argmin} \left\{ \sum_{j=1}^{B} [Y_{n+j} - \hat{m}_r(X_{n+j}; h)]^2 : h > 0 \right\}.$$

Then one can take $\hat{r}$ for a stopping rule defined by

$$\hat{r} = \operatorname{argmin}\left\{\sum_{j=1}^{B}[Y_{n+j} - \hat{m}_r(X_{n+j}; \hat{h}_r)]^2 : r \geq 1\right\}$$

and the data-driven bandwidth $\hat{h}_{\hat{r}}$. It would be of interest to see whether the regression estimator with these data-driven choices $\hat{r}$ and $\hat{h}_{\hat{r}}$ achieves the minimax optimal rate without $\nu$, the smoothness of the underlying function, being specified. We leave this as an open problem.

A method based on a cross-validatory criterion can be described similarly. As an alternative to these methods that are based on estimation of the prediction error, one may estimate the mean squared error of the estimator $\hat{m}_r(\cdot; h)$ and then choose $h$ and $r$ that minimize the estimated mean squared error. There have been many proposals for estimating the mean squared errors of kernel-based estimators of the regression function in connection with bandwidth selection; see, for example, Ruppert, Sheather and Wand (1995) and Section 4.3 of Fan and Gijbels (1996).

## 3. Numerical properties

In this section, we present the finite sample properties of the $L_2$ boosting estimates. To see how $L_2$ boosting compares favorably to the use of higher-order kernels as a method of bias reduction, we consider

$$\bar{m}_r(x) = \frac{\sum_{i=1}^{n} K_h^{[r]}(X_i - x) Y_i}{\sum_{i=1}^{n} K_h^{[r]}(X_i - x)},$$

where $K^{[r]}$ is a $2(r+1)$th-order kernel defined by, with $K^{[0]} = K$,

$$K^{[r]}(x) = 2K^{[r-1]}(x) - K^{[r-1]} * K^{[r-1]}.$$

Sufficient smoothness of $m$ and $f$ permitting, $\bar{m}_r$ is known to have a bias of order $h^{2(r+1)}$, which is of the same magnitude as the bias of the $r$-step boosted estimate $\hat{m}_r$.

The simulation was done under the following two models:

(1) $m(x) = \sin(2\pi x), \qquad 0 \leq x \leq 1;$

(2) $m(x) = \frac{2}{5}\{3\sin(4\pi x) + 2\sin(3\pi x)\}, \qquad 0 \leq x \leq 1.$

We took $U(0,1)$ for the distribution of $X_i$, and $N(0, 0.5^2)$ for the errors. For each model, two hundred pseudo-samples of size $n = 100$ and $400$ were generated. We used the Gaussian kernel $K$. We evaluated the mean integrated squared errors (MISE) of the estimates based on these samples. For this, we took 101 equally spaced grid points on [0,1] and used the trapezoidal rule for the numerical integration.

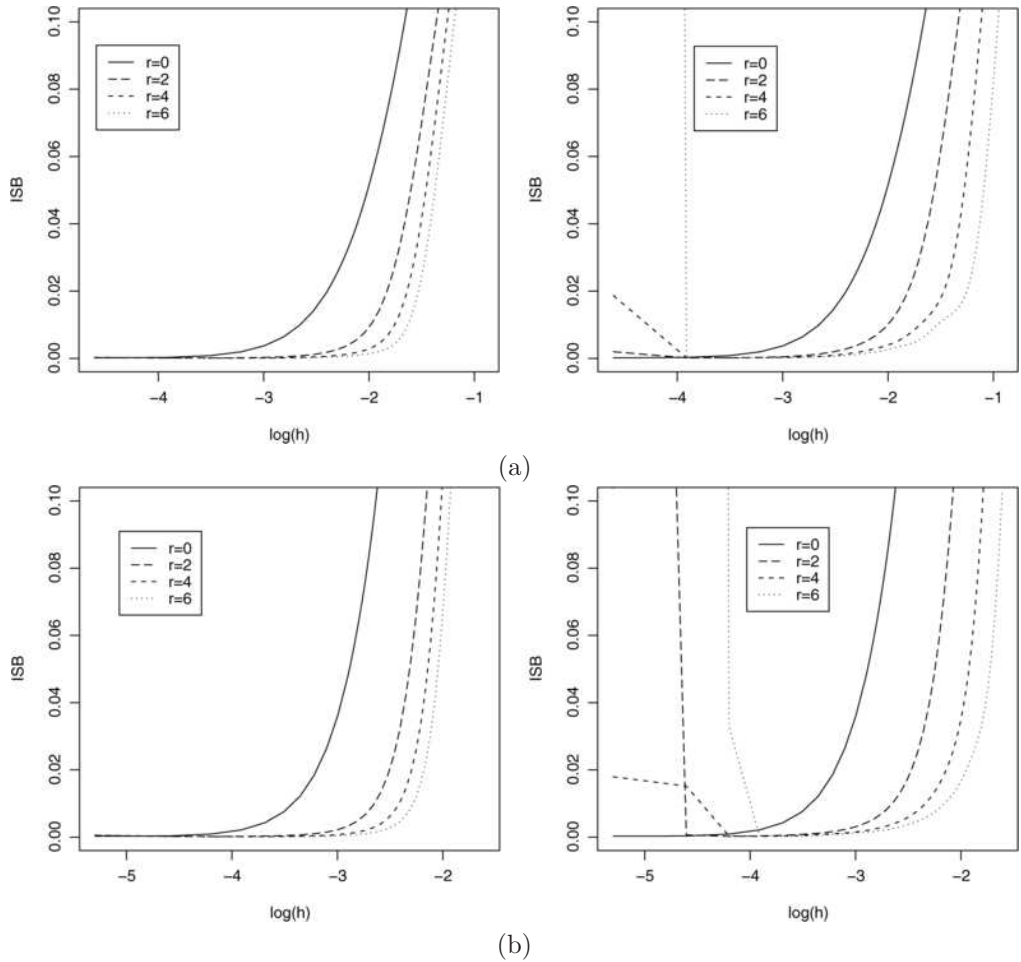

**Figure 1.** (a) Integrated squared bias for $r = 0, 1, \ldots, 6$ based on 200 pseudo-samples of size $n = 400$ from the model (1). (b) Integrated squared bias for $r = 0, 1, \ldots, 6$ based on 200 pseudo-samples of size $n = 400$ from the model (2). The left panel is for the $L_2$ boosting estimate and the right panel is for the higher-order kernel estimate.

Figures 1–3 show how the bias, variance and MISE of the estimates change as the boosting iteration number or the order of the kernel increases when $n = 400$. The result for $r = 0$ corresponds to the Nadaraya–Watson estimate. The curves in Figures 1 and 2 depict the integrated squared biases (ISB) and the integrated variance (IV), respectively, as functions of the bandwidth, and those in Figure 3 represent MISE. Table 1 gives the minimal MISE along with the optimal bandwidths that attain the minimal values for both sample sizes $n = 100$ and 400.

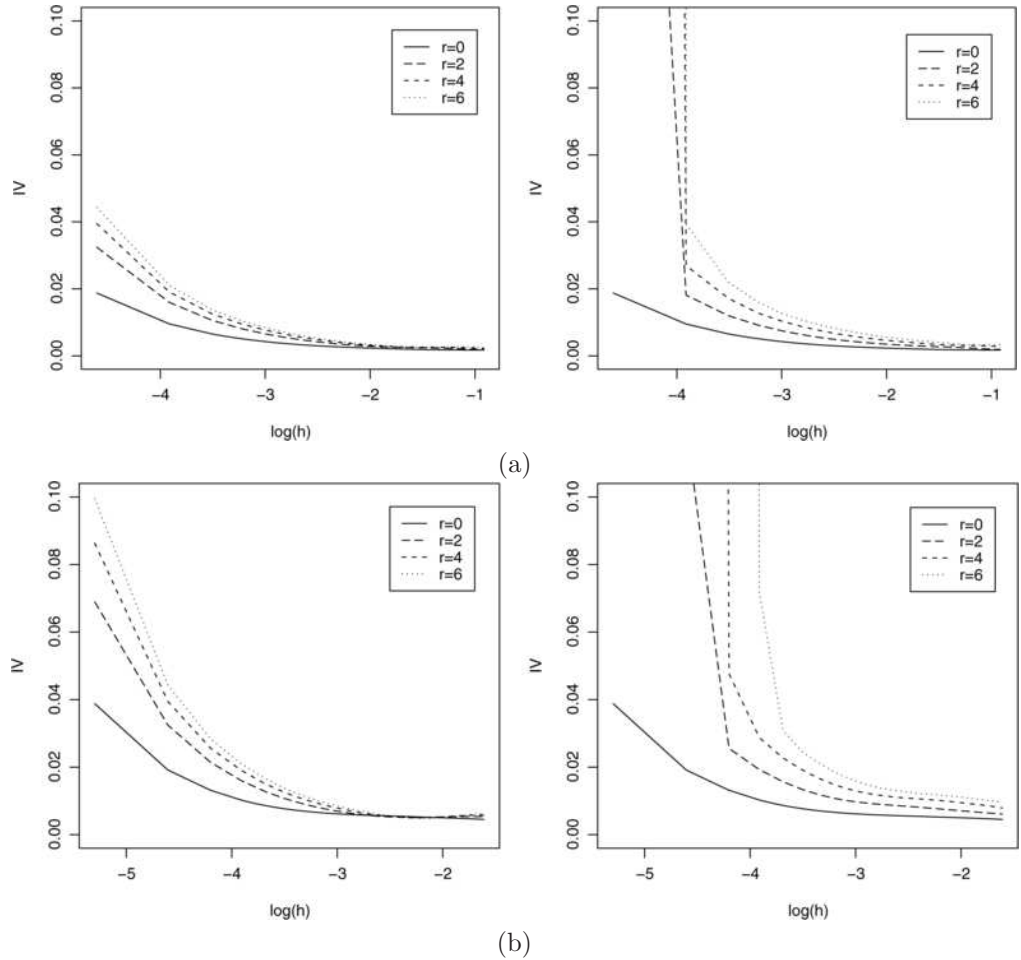

**Figure 2.** (a) Integrated variance for $r = 0, 1, \ldots, 6$ based on 200 pseudo-samples of size $n = 400$ from the model (1). (b) Integrated variance for $r = 0, 1, \ldots, 6$ based on 200 pseudo-samples of size $n = 400$ from the model (2). The left panel is for the $L_2$ boosting estimate and the right panel is for the higher-order kernel estimate.

For the $L_2$ boosted estimates, we see from the figures that the ISB reduces as the boosting iteration number $r$ increases in the whole range of the bandwidth. In particular, it decreases rapidly at the beginning of the boosting iteration and the degree of reduction decreases as $r$ increases. On the other hand, the IV increases at a relatively slower rate as $r$ increases. Since the decrement of the ISB (as $r$ increases) is greater than the increment of the IV for moderate-to-large bandwidths ($h \geq \mathrm{e}^{-3.0} \approx 0.05$ for model (1) and $h \geq \mathrm{e}^{-3.7} \approx 0.025$ for model (2)), and the former is smaller than the latter for small bandwidths, the value of MISE gets smaller as $r$ increases in the range of moderate-to-large bandwidths, while it becomes larger in the range of small bandwidths. The results in Table 1 show that the minimal value of MISE always decreases and the optimal bandwidth gets larger as $r$ increases. These results confirm our theoretical findings that $L_2$ boosting improves the order of the bias while not deteriorating the order of the variance.

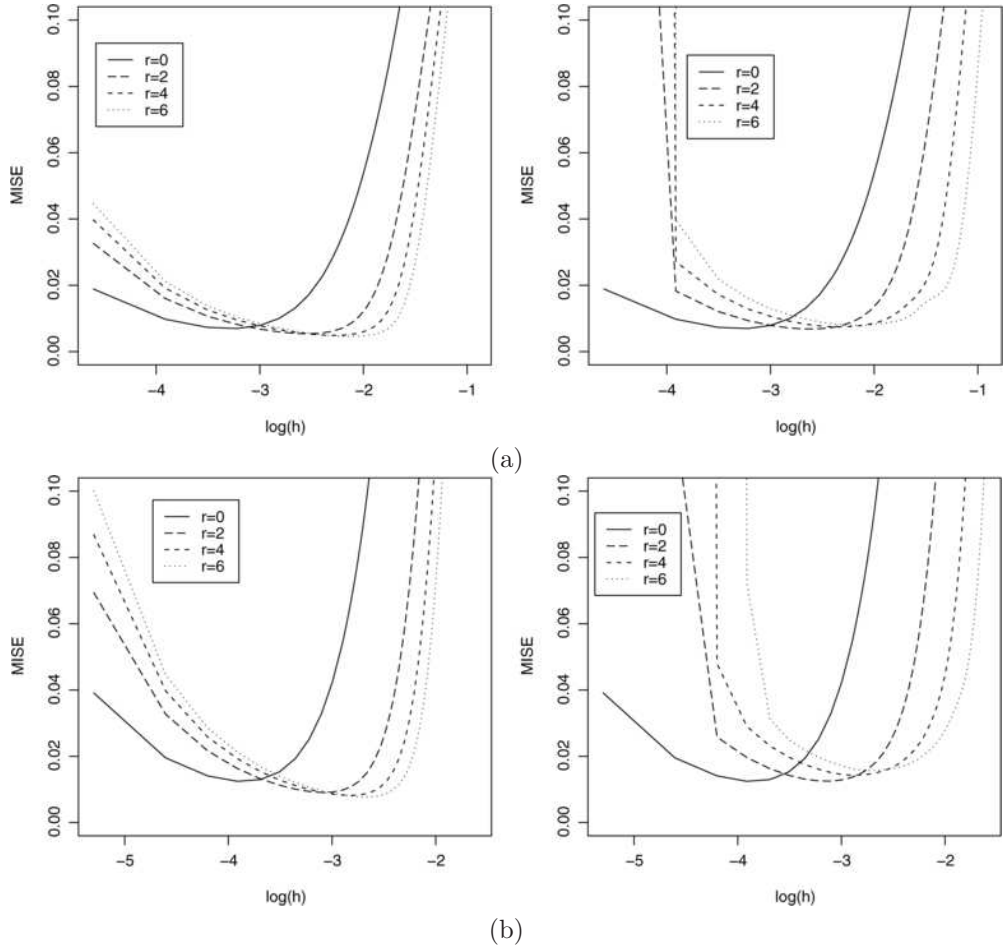

**Figure 3.** (a) Mean integrated squared error for $r = 0, 1, \ldots, 6$ based on 200 pseudo-samples of size $n = 400$ from the model (1). (b) Mean integrated squared error for $r = 0, 1, \ldots, 6$ based on 200 pseudo-samples of size $n = 400$ from the model (2). The left panel is for the $L_2$ boosting estimate and the right panel is for the higher-order kernel estimate.

For the higher-order kernel estimates, the behavior of the ISB and the IV as the order of kernel $r$ changes is similar to that of $L_2$ boosting except for small bandwidths. For small bandwidths, not only the IV but also the ISB increases as $r$ in-creases. This is contrary to the theory. In particular, the values of the ISB and IV explode when $r$ is large. Although not presented in this paper, we observed that the bad behavior is more severe when $n = 100$ and it starts at a relatively larger bandwidth than in the case of $n = 400$. Furthermore, Table 1 reveals that the minimal value of MISE starts to increase at some point as the order of the kernel $r$ increases. This erratic behavior of the higher-order kernel estimate is due to the fact that its denominator often takes near-zero or even negative values, which occurs more often for larger $r$, and it makes the estimate very unstable. This suggests that, contrary to $L_2$ boosting, the theoretical advantages of higher-order kernels do not take effect in practice.

**Table 1.** Minimal MISE and the corresponding optimal bandwidth $h$

| | Model (1) | | | | | Model (2) | | | |
| | $L_2$-boosting | | Higher-order kernel | | | $L_2$-boosting | | Higher-order kernel | |
| $r$ | $h$ | MISE | $h$ | MISE | $r$ | $h$ | MISE | $h$ | MISE |
| --- | --- | --- | --- | --- | --- | --- | --- | --- | --- |
| $n=100$ | | | | | | | | | |
| 0 | 0.050 | 0.0215 | 0.050 | 0.0215 | 0 | 0.030 | 0.0431 | 0.030 | 0.0431 |
| 1 | 0.080 | 0.0188 | 0.080 | 0.0208 | 1 | 0.045 | 0.0355 | 0.045 | 0.0436 |
| 2 | 0.100 | 0.0176 | 0.100 | 0.0213 | 2 | 0.060 | 0.0324 | 0.070 | 0.0493 |
| 3 | 0.120 | 0.0168 | 0.120 | 0.0223 | 3 | 0.065 | 0.0305 | 0.085 | 0.0544 |
| 4 | 0.130 | 0.0162 | 0.140 | 0.0231 | 4 | 0.075 | 0.0293 | 0.100 | 0.0588 |
| 5 | 0.140 | 0.0157 | 0.160 | 0.0238 | 5 | 0.080 | 0.0284 | 0.110 | 0.0624 |
| 6 | 0.150 | 0.0153 | 0.180 | 0.0242 | 6 | 0.085 | 0.0277 | 0.125 | 0.0649 |
| $n=400$ | | | | | | | | | |
| 0 | 0.040 | 0.0070 | 0.040 | 0.0070 | 0 | 0.020 | 0.0124 | 0.020 | 0.0124 |
| 1 | 0.060 | 0.0059 | 0.060 | 0.0066 | 1 | 0.035 | 0.0099 | 0.035 | 0.0118 |
| 2 | 0.080 | 0.0054 | 0.070 | 0.0068 | 2 | 0.045 | 0.0091 | 0.045 | 0.0125 |
| 3 | 0.090 | 0.0051 | 0.090 | 0.0072 | 3 | 0.055 | 0.0086 | 0.050 | 0.0134 |
| 4 | 0.100 | 0.0049 | 0.100 | 0.0075 | 4 | 0.060 | 0.0082 | 0.055 | 0.0143 |
| 5 | 0.110 | 0.0047 | 0.110 | 0.0077 | 5 | 0.065 | 0.0080 | 0.060 | 0.0150 |
| 6 | 0.120 | 0.0046 | 0.120 | 0.0079 | 6 | 0.070 | 0.0077 | 0.065 | 0.0157 |

# Acknowledgement

This work was supported by the Korea Research Foundation Grant funded by the Korean Government (MOEHRD)(KRF-2005-070-C00021).